\documentclass[a4paper, oneside, 11pt]{article}
\usepackage[english]{babel}
\usepackage{latexsym,amsfonts,amsmath,enumerate,graphics,enumerate,amsthm}
\title{\textbf{Topological concepts in partially ordered vector spaces} }
\author{T.\ Hauser}

\let\epsilon=\varepsilon

\theoremstyle{definition}
\newtheorem{definition}{Definition}[section]
\newtheorem{theorem}[definition]{Theorem}
\newtheorem{proposition}[definition]{Proposition}
\newtheorem{lemma}[definition]{Lemma}
\newtheorem{corollary}[definition]{Corollary}
\theoremstyle{remark}
\newtheorem{remark}[definition]{Remark}
\newtheorem{example}[definition]{Example}

\newtheorem*{acknowledgement}{Acknowledgement}

\let\epsilon=\varepsilon
\let\phi=\varphi
\let\theta=\vartheta
\newcommand{\R}{\mathbb{R}}
\newcommand{\N}{\mathbb{N}}

\begin{document}
\maketitle

\begin{abstract}
	In the context of partially ordered vector spaces one encounters different sorts of order convergence and order topologies. This article will investigate these notions and their relations. In particular we study and relate the order topology presented by Floyd, Vulikh and Dobbertin, the order bound topology studied by Namioka and the concept of order convergence given in the works of Abramovich, Sirotkin, Wolk and Vulikh. 
\end{abstract}

\section{Introduction}

In the study of partially ordered vector spaces one uses topological concepts like order convergence and order continuity as can be seen for example in   \cite{Abra,Wolk,Vulikh67,Positiveoperators,CAD}. In particular one encounters different types of order convergence, which lead to different types of order continuity as well as different types of order topology. The different types of order continuity are for example studied in \cite{Abra} and \cite{Paper1} and we will focus in this paper on properties and relationships of different sorts of topologies that can be defined in partially ordered vector spaces. In particular we will investigate properties and relationships of the \emph{order topology} \cite{Floyd1955,Vulikh67,Dobbertin84,Imh} and the \emph{order bound topology} \cite{Namioka57,CAD}. We start by defining the order topology of a partially ordered vector space $X$. Recall\footnote{Missing definitions and further comments can be found in the preliminaries below.} that we write $x_\alpha\downarrow 0$, if $\inf_{\alpha\in A}x_\alpha=0$ and if $x_\alpha$ is decreasing, i.e.\ $\alpha\mapsto -x_\alpha$ is monotone. 

\begin{definition}
	$M\subseteq X$ is called \emph{net catching for $x\in X$}, if for all nets $(\check{x}_\alpha)_{\alpha \in A}$ with $\check{x}_\alpha \downarrow 0$ there is $\alpha \in A$ such that $[x-x_\alpha,x+x_\alpha]\subseteq M$. 
	A subset $O$ of $X$ is called \emph{order open} if $O$ is a net catching set for every $x\in O$. The \emph{order topology} is defined as $\tau_o(X):=\{O\subseteq X;\, O \mbox{ is order open}\}$.
\end{definition}

	It is natural to ask whether there are always non empty, order open and order bounded sets in a partially ordered vector space. 
	In Section 2 we will see that this property is equivalent to the existence of order topological interior points of the cone $X_+$. We will thus characterise the interior of the cone with respect to the order topology as the set of all net catching elements, a notion defined as follows. 
	An element $y \in P$ is called \emph{net catching}, if for each net $(\check{x}_\alpha)_{\alpha \in A}$ with $\check{x}_\alpha \downarrow 0$ there is some $\alpha \in A$ such that $\check{x}_\alpha\leq y$.
	Section 3 is devoted to the interplay of net catching elements and order units of $X$. It is shown that whenever $X$ contains net catching elements, then every order unit is net catching. Furthermore we present that whenever $X$ is Archimedean and directed, then every net catching elements is an order unit. This section ends with examples.

	We then ask for necessary conditions on a partially ordered vector space $X$ under which one obtains the existence of net catching elements, i.e.\ non empty, order open and order bounded sets. In Section 4 we present such a condition. In particular we show in Theorem \ref{the:necessarryconditionsfornetcatchingelements} that whenever there is a linear topology $\tau$ on $X$ and a $\tau$-compact base $B$ of $X_+$, then each upper bound of $B$ is net catching. 
	This condition will give us the tool to establish the existence of net catching elements for several examples. In Section 5 we present that Archimedean finite dimensional partially ordered vector spaces and reflexive Banach spaces equipped with ice cream cones have net catching elements. 
	In particular we obtain all finite dimensional Archimedean vector lattices to contain non empty, order open and order bounded sets. In Section 6 we will ask whether this is also true for all infinite dimensional Archimedean vector lattices. We show in that section that the opposite is valid and present in Theorem \ref{the:Characterization_existencence_for_Archimedean_VL} that an Archimedean vector lattice is finite dimensional if and only if there exists a non empty, order open and order bounded set in $X$. 
	
	 In \cite{Namioka57} and \cite{CAD} the concept of the \emph{order bound topology} is introduced and studied. We present in Section 7 the definition of this concept and that it differs from the concept of order topology in general. We investigate the relationship of these topologies and the following concept of order convergence.
Let $x \in X$ and let $(x_\alpha)_{\alpha \in A}$ be a net in $X$. We define $x_\alpha \xrightarrow{o} x$, 
		if and only if there is a net  $(\check{x}_\beta)_{\beta \in B}$
		and 
		a mapping 
		$\eta\colon B \rightarrow A$ such that $\check{x}_\beta \downarrow 0$ and
		$\pm (x_\alpha-x) \leq \check{x}_\beta$ for every 
		$\beta\in B$ and $\alpha \in A$ with $\alpha\geq \eta(\beta)$.
Note that this concept is investigated under the notion of $o_3$-convergence in \cite{Paper1} and also studied in \cite{Abra,Wolk,Vulikh67,Imh}. 	 
	 We furthermore show that the order units play an analogous role for the order bound topology as the net catching elements for the order topology.
	 In Section 8 we will see that under the assumption of the existence of net catching elements in $X$ the order topology $\tau_o$ is normable, if and only if $X$ is Archimedean and directed, and if and only if $\tau_o$ equals the order bound topology. Furthermore we will obtain under this assumption that the order topological convergence is equivalent to the concept of order convergence defined above. 
	 We use these insights in Section 9 to present that the topology of a reflexive Banach space can be obtained from the order topology of any ice cream cone. 
	In the last section we will consider the following problem. Whenever a net 
	is decreasing to $0$, then we know that this net also converges to $0$ with respect to the order topology. It is thus natural to expect that an analogue is true for the order bound topology. Nevertheless we will see in Theorem \ref{the:characterizationwhendodecreasingnetsconvergewithrespecttotheorderboundtopology} that in Archimedean partially ordered vector spaces that contain order units this property is equivalent to the existence of net catching elements, i.e.\ to the existence of non empty order open and order bounded sets. Thus in particular we can observe that in infinite dimensional Archimedean vector lattices with order units there is always a net that decreases to $0$ but which does not converge with respect to the order bound topologies. 
	
We next fix some notation. As usual, on a non-empty set $P$ a binary relation $\leq$ is called a \emph{partial order} if it is reflexive, transitive and anti-symmetric. The set $P$ is then  called a \emph{partially ordered set}. 
For $x,y\in P$ we write $x<y$ if $x\leq y$ and $x\neq y$.
For $U,V\subseteq P$ we write $U\leq V$ if for every $u\in U$ and $v\in V$ we have $u\leq v$. If $V=\{v\}$ for $v\in P$, we abbreviate $U\leq \{v\}$ by $U\leq v$. For $u\in P$, $u\leq V$ is defined similarly.
For $x\in P$ we denote $P_{\geq x}:=\{z\in P;\ x\leq z\}$.
For $x,y\in P$ the \emph{order interval} is given by $[x,y]:=\{z\in P; \, x\leq z\leq y\}$. 
A set $M\subseteq P$ is called \emph{full} if for every $x,y\in M $ one has $[x,y]\subseteq M$.
$P$ is called \emph{directed (upward)} if for every $x,y\in P$ the set $P_{\geq x}\cap P_{\geq y}$ is non-empty. \emph{Directed downward} is defined analogously. 
For a subset of $P$,  the notions \emph{bounded above}, \emph{bounded below},
\emph{order bounded},
  \emph{upper (or lower) bound} and  \emph{infimum (or supremum)} 
are defined as usual.

We call $P$ a \emph{lattice} if for every non-empty finite subset of $P$ the infimum and the supremum exist in $P$, and a \emph{Dedekind complete lattice} if  
every non-empty set that is bounded above has a supremum, and every non-empty set that is bounded below has an infimum.  

A set $M\subseteq P$ is called \emph{majorizing} in $P$ if for every $y\in P$ the set $P_{\geq y}\cap M$ is non-empty, and
\emph{order dense} in $P$ if for every $y\in P$ one has \[\sup\{x\in M;\, x\leq y\}=y=\inf\{x\in M;\, y\leq x\}.\]
 

Let $P$ and $Q$ be partially ordered sets and $f\colon P\to Q$ a mapping. $f$ is called \emph{monotone} if for every $x,y\in P$ with $x\leq y$ one has that $f(x)\leq f(y)$, and 
\emph{order reflecting} if for every $x,y\in P$ with
 $f(x)\leq f(y)$ one has that $x\leq y$.
Note that every order reflecting mapping is injective.
We call $f$ an \emph{order embedding} if $f$ is monotone and order reflecting.
The map $f$ is called an \emph{order isomorphism}, whenever $f$ is bijective, monotone and order reflecting. 
For sets $U\subseteq P$ and $V\subseteq Q$ we use the notation $f[U]$ for the image of $U$ under $f$, and $[V]f$ for the preimage of $V$ under $f$.

On a real vector space $X$, we consider a partial order $\leq$ on $X$ such that for every $\lambda \in \mathbb{R}_+$, $x,y,z\in X$ with $x\leq y$ one has that $x+z\leq y+z$ and $\lambda x\leq \lambda y$. Then $X$ is called a \emph{partially ordered vector space}. We denote $X_+:=X_{\geq 0}$ and call the elements of $X_+$ \emph{positive}. $X_+$ is called \emph{generating} if $X=X_+-X_+$. Note that $X$ is directed if and only if $X_+$ is generating. 
A non empty convex subset $B$ of $X_+\setminus\{0\}$ is called a \emph{base (of $X_+$)}, if for each $x\in X_+\setminus\{0\}$ there are unique $b\in B$ and $\lambda\in \Lambda$ such that $x=\lambda b$. 
We say that $X$ is \emph{Archimedean} if for every $x,y\in X$ with $nx\leq y$
for all $n\in \mathbb{N}$ one has that $x\leq 0$. An element $u\in X_+$ is called an \emph{order unit}, if for all $x\in X$ there is $n\in \mathbb{N}$ such that $x\leq n u$. If $X$ is a lattice, we call $X$ a \emph{vector lattice}.
	 A linear mapping $f\colon X \to Y$ between partially ordered vector spaces is said to be \emph{positive}, if it is monotone. 

If $X$ is an Archimedean and directed partially ordered vector space, then there is a Dedekind complete vector lattice $X^\delta$ and a linear order embedding $J\colon X \to X^\delta$ such that $J[X]$ is order dense in $X^\delta$. As usual, $X^\delta$ is called \emph{Dedekind completion} of $X$.



Next we define the order topology of a partially ordered set $P$. Note that this definition and the definition of net catching sets are equivalent to the definition above for partially ordered vector spaces as presented in \cite[Prop. 5.1]{Paper1}. 
A subset $M\subseteq P$ is called a 
	\emph{net catching set}
	for $x\in P$ if for all nets  $(\hat{x}_\alpha)_{\alpha\in A}$ and $(\check{x}_\alpha)_{\alpha\in A}$ in $P$
	with $\hat{x}_\alpha \uparrow x$ and $\check{x}_\alpha\downarrow x$ there is
	$\alpha \in A$ such that $[\hat{x}_\alpha,\check{x}_\alpha]\subseteq M$.
	A subset $O$ of $P$ is called \emph{order open} if $O$ is a net catching set for every $x\in O$. A subset $C$ of $P$ is called \emph{order closed} if $P\setminus C$ is order open. The \emph{order topology} is defined as $\tau_o(P):=\{O\subseteq P;\, O \mbox{ is order open}\}$.

\begin{remark} \label{rem:basicsonordertopology}
\begin{itemize}
\item[(a)] For every  $p\in P$ we know $P_{\geq p}$ to be order closed. Hence the cone $X_+$ is order closed \cite[Corollary 3.19]{Paper1}.
	\item[(b)]  Let $(x_\alpha)_{\alpha\in A}$ be a net in $P$ with $x_\alpha \xrightarrow{\tau_o} x\in P$. If $x_\alpha \downarrow$, then also $x_\alpha \downarrow x$. This follows from \cite[Lemma 4.3]{Paper1}.
	\item[(c)] A subset $M\subseteq X$ is net catching for $x\in X$, if and only if $M-x$ is net catching for $0$. 
\end{itemize} 
\end{remark}

Also $o$-convergence can be defined and studied in partially ordered sets. 
Let $x \in P$ and let $(x_\alpha)_{\alpha \in A}$ be a net in $P$. We define
	$x_\alpha \xrightarrow{o} x$, 
	if there are nets $(\hat{x}_\beta)_{\beta \in B}$
	and $(\check{x}_\gamma)_{\gamma \in C}$ in $P$ and a mapping 
	$\eta\colon B \times C \rightarrow A$ such that $\hat{x}_\beta \uparrow x$,
	$\check{x}_\gamma \downarrow x$ and $\hat{x}_\beta \leq x_\alpha \leq \check{x}_\gamma$ for every 
	$\beta\in B$, $\gamma \in C$ and $\alpha \in A_{\geq \eta(\beta,\gamma)}$. Note that this concept is the concept of $o_3$-convergence in \cite{Paper1}.
	
\begin{remark} \label{rem:basicorderconvergencestuff}
	(a) Let $i\in \{1,2,3\}$ and $C\subseteq P$. Then $C$ is order closed if and only if for every net $(x_\alpha)_{\alpha \in A}$ with $x_\alpha \xrightarrow{o} x\in P$ it follows that $x\in C$. This is shown in \cite[Theorem 3.14]{Paper1}.
	
	(b)	Let	$(x_\alpha)_{\alpha \in A}$ be a net in $P$ and  $x \in P$. The 	following results can be found in \cite[Proposition 3.6]{Paper1} and \cite[Remark 3.7]{Paper1}.
		\begin{itemize}
			\item[(i)] $x_\alpha \downarrow x$ implies $x_\alpha \xrightarrow{o} x$,
			\item[(ii)] $x_\alpha \xrightarrow{o} x$ implies $x_\alpha \xrightarrow{\tau_o} x$.
	\end{itemize}

\end{remark}
 
A map $f\colon P\to Q$ between partially ordered sets $P$ and $Q$ is called \emph{$o$-continuous}, if for every net $(x_\alpha)_{\alpha\in A}$ and every $x\in P$ with $x_\alpha\xrightarrow{o}x$ we have that $f(x_\alpha)\xrightarrow{o}f(x)$. Furthermore $f\colon P \to Q$ is called \emph{order continuous}, if it is continuous with respect to the order topologies $\tau_o(P)$ and $\tau_o(Q)$, respectively.	
	
\begin{remark} \label{rem:characterisationordercontinuityforpositivemaps}
	(a) Let 
	  $f\colon P \to Q$ be a monotone map. Then $f$ is $o$-continuous if and only if $f$ is order continuous. This result can be found in \cite[Theorem 4.4]{Paper1}. If $f$ is an order embedding such that $f[P]$ is order dense in $Q$, then $f$ is order continuous and for every net $(x_\alpha)_{\alpha \in A}$ and every $x \in P$ we have $x_\alpha \downarrow x$ if and only if $f(x_\alpha)\downarrow f(x)$. These statements are proven in \cite[Corollary 4.5]{Paper1} and \cite[Proposition 1.3(ii)]{Paper1}.
	
	(b) Let $X$ and $Y$ be a partially ordered vector spaces and $f\colon X \to Y$ be a linear mapping. Then $f$ is positive and order continuous if and only if for every net $(x_\alpha)_{\alpha\in A}$ with $x_\alpha\downarrow 0$ it holds $f(x_\alpha)\downarrow 0$. This is shown in \cite[Proposition 7.1]{Paper1}. 
\end{remark}

\section{Net catching elements and the interior of the cone}

In this section let $P$ be a partially ordered set and $X$ a partially ordered vector space. We characterise the interior of $X_+$ with respect to the order topology. This question is motivated by the following proposition. 

\begin{proposition} \label{pro:characterisation_existenceorderopenandorderbddsets}
	There exists a non empty, order bounded and order open set in $X$, if and only if the order topologically interior $\tau_o\operatorname{-int}(X_+)$ of $X_+$ is not empty. 
\end{proposition}

\begin{proof}
	For $x \in \tau_o\operatorname{-int}(X_+)=:O$ let $U:=((2x)-O)\cap O$ to obtain 
	$x \in U\in \tau_o(X)$ and $U\subseteq [0,2x]$. To show the converse implication assume $O \in \tau_o(X)$ to be order bounded and non-empty. Thus there are $x,y,z \in X$ such that $z \in O \subseteq [x,y]$, hence $z-x \in O-x\subseteq [0,y-x] \subseteq X_+$. Since $O-x\in  \tau_o(X)$ this proves $z-x\in \tau_o\operatorname{-int}(X_+)$.
\end{proof}

\begin{definition}
	An element $y \in P$ is called \emph{downward net catching} with respect to $x\in P$, if for each net $(\check{x}_\alpha)_{\alpha \in A}$ with $\check{x}_\alpha \downarrow x$ there is some $\alpha \in A$ such that $\check{x}_\alpha\leq y$. Analogously we define \emph{upward net catching} with respect to $x\in P$.	
	Note that with the definition given in the introduction an element $x\in X$ is net catching, if and only if it is downward net catching with respect to $0$.	
\end{definition}

\begin{remark} \label{rem:ncelements}
	(i) Let $x,y,z \in P$. If $y$ is downward net catching with respect to $x$ and $z\geq y$, then $z$ is downward net catching with respect to $x$. Thus if $x\in G$ is net catching and $y\geq x$, then $y$ is net catching as well. 
	
	(ii) Let $x,y,z \in X$. Clearly if $x$ is net catching with respect to $y$, then $x+z$ is net catching with respect to $y+z$. Note furthermore that $x$ is upward net catching with respect to $y$ if and only if $y$ is downward net catching with respect to $x$. Indeed, if $x$ is upward net catching with respect to $y$ and if $\check{x}_\alpha \downarrow x$, then $-\check{x}_\alpha+x+y \uparrow y$ and so there is $\alpha \in A$ such that $-\check{x}_\alpha+x+y\geq x$, i.e.\ $\check{x}_\alpha\leq y$. This proves $y$ to be downward net catching with respect to $x$. 
	
	(iii) Let $\lambda \in R_{\geq 0}\setminus \{0\}$.\\
	(iiia) If $x\in X$ is downward net catching with respect to $y\in X$, then $\lambda x$ is downward net catching with respect to $\lambda y$. Indeed, if $(\check{x}_\alpha)_{\alpha\in A}$ is a net in $X$ with $\check{x}_\alpha \downarrow \lambda y$, then $\frac{1}{\lambda}\check{x}_\alpha \downarrow y$. Thus there is $\alpha \in A$ such that $\frac{1}{\lambda}\check{x}_\alpha \leq x$, i.e.\ $\check{x}_\alpha \leq \lambda x$. A similar argument shows the analogue if $x$ is upward net catching with respect to $y$. \\
	(iiib) If $x\in X$ is net catching, so is $\lambda x$. 
\end{remark}

The relationship between downward and upward net catching elements and net catching sets is given by the following Proposition. 

\begin{proposition} \label{pro:nc_elements_and_nc_sets}
	If $x,y,z\in P$, then $[x,z]$ is a net catching set for $y$ if and only if $x$ is upward net catching and $z$ is a downward net catching with respect to $y$. If $x\in X$, then $x$ is net catching if and only if $[-x,x]$ is a net catching set for $0$. 
\end{proposition}

\begin{proof}
	Assume that $[x,z]$ is a net catching set for $y$. Let $(\check{x}_\alpha)_{\alpha \in A}$ be a net in $P$ with $\check{x}_\alpha \downarrow y$. Setting $\hat{x}_\alpha:=y$ for all $\alpha \in A$ yields a net $(\hat{x}_\alpha)_{\alpha \in A}$ with $\hat{x}_\alpha \uparrow y$. Thus since $[x,z]$ is a net catching set for $y$ there is $\alpha \in A$ such that $[\hat{x}_\alpha,\check{x}_\alpha]\subseteq [x,z]$, hence $\check{x}_\alpha\leq z$. This proves $z$ to be downward net catching with respect to $y$. A similar proof shows $x$ to be upward net catching with respect to $y$. 
	
	To show the converse implication assume $x$ to be upward net catching and $z$ to be a downward net catching with respect to $y$. Furthermore let $(\hat{x}_\alpha)_{\alpha \in A}$ and $(\check{x}_\alpha)_{\alpha \in A}$ be nets in $P$ with $\hat{x}_\alpha \uparrow y$ and $\check{x}_\alpha \downarrow y$. Since $x$ is upward net catching and $z$ is downward net catching with respect to $y$ there are $\alpha_x,\alpha_z\in A$ such that $x\leq \hat{x}_{\alpha_x}$ and $\check{x}_{\alpha_z}\leq z$. Let $\alpha \in A_{\geq \alpha_x}\cap A_{\geq \alpha_z}$ to obtain $x\leq \hat{x}_{\alpha_x}\leq \hat{x}_\alpha \leq \check{x}_\alpha \leq \check{x}_{\alpha_z}\leq z$, i.e.\ $[\hat{x}_\alpha,\check{x}_\alpha]\subseteq [x,z]$.
	
	The second statement is a direct consequence of the first one, since in a partially ordered vector space an element $x$ is downward net catching with respect to $0$ if and only if $-x$ is upward net catching with respect to $0$.  
\end{proof}

\begin{lemma}\label{lem:interior_pts_are_ptsst_p_is_upward_nc}
	Let $p \in P$ and $x$ be an order topological interior point of $P_{\geq p}$. Then $p$ is upward net catching with respect to $x$.
\end{lemma}

\begin{proof}
	Let $(\hat{x}_\alpha)_{\alpha \in A}$ be a net in $P$ with $\hat{x}_\alpha \uparrow x$. Setting $\check{x}_\alpha:=x$ we obtain a net $(\check{x}_\alpha)_{\alpha \in A}$ with $\check{x}_\alpha \downarrow x$. Since $\tau_o\operatorname{-int}(P_{\geq p})$ is a net catching set for $x$ we are now able to conclude the existence of $\alpha \in A$ such that $[\hat{x}_\alpha,\check{x}_\alpha]\subseteq \tau_o\operatorname{-int}(P_{\geq p})\subseteq P_{\geq p}$, hence $p\leq \hat{x}_\alpha$. 
\end{proof}

	Combining Lemma \ref{lem:interior_pts_are_ptsst_p_is_upward_nc} with Remark \ref{rem:ncelements}(ii) we obtain that in $X$ every interior point of the cone $X_+$ is net catching. We will see in the next theorem that also the converse holds true. 
	

\begin{theorem} \label{the:interiorofconeare_nc_elements}
	Let $X$ be a partially ordered vector space. Then
	\[\tau_o\operatorname{-int}(X_+)=\{x \in X;\, x \text{ net catching}\}.\] 
\end{theorem}

\begin{proof}
	 In order to show the remaining inclusion it is sufficient to show that $U:=\{x \in X;\, x \text{ net catching}\}$ is order open. Let $x\in X$ be net catching and set $y:=\frac{1}{2}x$. Then $y$ is net catching as well and from Proposition \ref{pro:nc_elements_and_nc_sets} we obtain $[-y,y]$ to be a net catching set for $0$. Thus $[x-y,x+y]$ is a net catching set for $x$. By Remark \ref{rem:ncelements}(i) we conclude $[x-y,x+y]=[y,x+y]\subseteq U$, which shows $U$ to be a net catching set for $x$. This proofs $U\in \tau_o(G)$. 
\end{proof}

%

\section{Net catching elements and order units}

In this section let $X$ be a partially ordered vector space. This section is devoted to the presentation of the interplay of net catching elements and order units of $X$.

\begin{proposition} \label{pro:ncexist_thenorderunitsare_nc}
	If $X$ contains net catching elements, then every order unit in $X$ is net catching. 
\end{proposition}

\begin{proof}
	Let $u$ be an order unit of $X$ and $x\in X$. There is $n \in \mathbb{N}$ such that $x\leq n u$. If $x$ is net catching, then $n u$ is net catching as well by Remark \ref{rem:ncelements}(i).  
	The statement follows from Remark \ref{rem:ncelements}(iiib).
\end{proof}

\begin{remark}
	In Example \ref{exa:l_inftyhasnonetcatchingelements} we introduce a Dedekind complete vector lattice with order units but no net catching elements. Thus the assumption on the existence of net catching elements can not be omitted in Proposition \ref{pro:ncexist_thenorderunitsare_nc}.  

	An analogue of the second statement in Proposition \ref{pro:ncexist_thenorderunitsare_nc} with the interchanged roles of order units and net catching elements is not valid in general. In Example \ref{exa:lexicographicsequencespacenetcatchingelements} we present a vector lattice with order units that contains net catching elements which are no order units.
\end{remark}

\begin{example} \label{exa:l_inftyhasnonetcatchingelements}
	Let $l_\infty$ be the vector space of all bounded and real sequences, equipped with the component wise order. Note that $l_\infty$ is a Dedekind complete vector lattice and $e:=(1)_{k \in \mathbb{N}}$ is an order unit of $l_\infty$. Furthermore $l_\infty$ has no net catching elements. Indeed, assume $l_\infty$ to have net catching elements. Proposition \ref{pro:ncexist_thenorderunitsare_nc} implies $e$ to be net catching. For $n,k\in \mathbb{N}$ let $e^{(n)}_k:=0$, if $k<n$ and $e^{(n)}_k:=2$ otherwise. This defined a sequence $(e^{(n)})_{n \in \mathbb{N}}$ in $l_\infty$ such that $e^{(n)}\downarrow 0$. Since $e^{(n)}\leq e$ for no $n \in \mathbb{N}$ and $e$ is assumed to be net catching we obtain a contradiction. 
\end{example}

\begin{proposition} \label{pro:nce_are_orderunits_XArchm_directed}
	Assume $X$ to be directed and Archimedean. Then each net catching element of $X$ is an order unit. 
\end{proposition}

\begin{proof}
	Let $x\in X$ be net catching and $y\in X$. As $X$ is directed there exists $z\in X_+$ such that $z\geq y$. The Archimedean property of $X$ implies $\frac{1}{n}z \downarrow 0$. Thus since $X$ is net catching there is $n \in \mathbb{N}$ such that $\frac{1}{n}z\leq x$, i.e.\ $y\leq z \leq n x$. 
\end{proof}

\begin{remark}
	(i) If $X$ is an Archimedean and directed partially ordered vector space we can characterise the set of all order topological interior points of the cone as the set of all net catching order units. This follows from Theorem \ref{the:interiorofconeare_nc_elements} and Proposition \ref{pro:nce_are_orderunits_XArchm_directed}.
	
	(ii) Note that Example \ref{exa:l_inftyhasnonetcatchingelements} provides a directed and Archimedean partially ordered vector space and an order unit that is not net catching. Furthermore in Proposition \ref{pro:finitedimensionalandArchimedeanspaceshavenetcatchinglements} we present that each Archimedean finite dimensional partially ordered vector space has net catching elements. Thus the directedness of $X$ can not be omitted in Proposition \ref{pro:nce_are_orderunits_XArchm_directed}. Moreover the following Example provides a vector lattice with net catching elements that are no order units. Thus the Archimedean property of $X$ can not be omitted in Proposition \ref{pro:nce_are_orderunits_XArchm_directed} as well.
	
\end{remark}

\begin{example} \label{exa:lexicographicsequencespacenetcatchingelements}
	Let $X$ be a non trivial vector subspace of the space of all sequences $(x_n)_{n \in \mathbb{N}}$ in $\R$ and ordered lexicographic, i.e.\ $(x_n)_{n \in \mathbb{N}} \leq (y_n)_{n \in \mathbb{N}}$, if $(x_n)_{n \in \mathbb{N}}=(y_n)_{n \in \mathbb{N}}$ or if there is $n \in \N$ such that $x_m=y_m$ for all $m<n$ and $x_n < y_n$. One easily verifies, that $\geq$ is a total vector ordering on $X$. Thus $X$ equipped with $\leq$ is a vector lattice. Note that
		\begin{align*}
			X_+:=\{0\} \cup \bigcup_{n \in \N}\{(x_n)_{n \in \mathbb{N}} \in X;\,0=x_m \text{ for all } m < n \text{ and } 0 < x_n \}.
		\end{align*}
We will show that $X_+\setminus \{0\}$ is the set of net catching elements and that $\{(x_n)_{n \in \mathbb{N}}\in X_+;\, 0<x_1\}$ is the set of order units of $X$. Note that $X$ can not be Archimedean by Proposition \ref{pro:nce_are_orderunits_XArchm_directed}.

Indeed, let $x\in X_+\setminus \{0\}$. To show that $x$ is net catching let $(x_\alpha)_{\alpha \in A}$ be a net in $X$ with $x_\alpha \downarrow 0$. If there would be no $\alpha \in A$ with $x_\alpha \leq x$, then $x < x_\alpha$ for all $\alpha \in A$, since $\leq$ is total. Thus $0 < x \leq \inf_{\alpha \in A}x_\alpha =0$, a contradiction. Hence there must be $\alpha \in A$ with $x_\alpha \leq x$.

To show the second statement let $x:=(x_n)_{n \in \mathbb{N}}\in X_+$ such that $0<x_1$. Thus for $y:=(y_n)_{n \in \mathbb{N}}\in X$ there is $m\in \mathbb{N}$ such that $y_1<m x_1$ and we have shown $y\leq m x$. If $x:=(x_n)_{n \in \mathbb{N}}$ is an order unit of $X$, then there is $m \in \mathbb{N}$ such that $(1)_{n \in \mathbb{N}}\leq m x$. Hence either $1=m x_1$ or $1<m x_1$ and we  obtain $0<x_1$. 
\end{example}

\section{Necessary conditions for net catching elements}
	
	In this section let $X$ be a partially ordered vector space. We will present that whenever there is a linear topology $\tau$ on $X$ and a $\tau$-compact base $B$ of $X_+$, then each upper bound of $B$ is net catching. 
	Recall\footnote{For a proof see \cite[Theorem 1.47]{CAD}} that a cone $X_+$ has a base if and only if $X$ admits a strictly positive linear functional $f\colon X\to \mathbb{R}$. If $f\colon X \to \mathbb{R}$ is any strictly positive functional, then $B_f=\{x \in X_+;\, f(x)=1\}$ is a base for $X_+$. Moreover if $B$ is a base of $X_+$, then there is a unique strictly positive functional $f_B\colon X \to \mathbb{R}$ such that $B=B_{f_B}$. We will now see, that a base $B$ is closed with respect to $\tau_o$ if and only if the corresponding strictly positive functional $f_B$ is order continuous. 
	
	\begin{proposition} \label{pro:Borderclosedif_f_ordercontinuous}
		Let $B$ a base of $X_+$ and $f\colon X \to \mathbb{R}$ be a strictly positive linear functional such that $B=\{x\in X_+;\, f(x)=1\}$. Then $B$ is order closed if and only if $f$ is order continuous.  
	\end{proposition} 
	
	\begin{proof}
		Note that $X_+$ is order closed by Remark \ref{rem:basicsonordertopology}(a). Thus if $f$ is order continuous, then $B=X_+ \cap f^{-1}(\{1\})$ is order closed. To show the converse implication assume $B$ to be order closed. We will use Remark \ref{rem:characterisationordercontinuityforpositivemaps}(b) to show, that $f$ is order continuous. Let $(x_\alpha)_{\alpha \in A}$ be a net in $X$ with $x_\alpha \downarrow 0$. Clearly $(f(x_\alpha))_{\alpha \in A}$ is a net in $\mathbb{R}_+$ such that $f(x_\alpha) \downarrow$. Let $\lambda$ be a lower bound of $\{f(x_\alpha);\, \alpha \in A\}$. To get a contradiction assume $\lambda >0$. Thus there is $f(x_\alpha)\geq \lambda >0$ for every $\alpha \in A$ and we define $y_\alpha:= \frac{1}{f(x_\alpha)}x_\alpha$ for all $\alpha \in A$ to obtain a net $(y_\alpha)_{\alpha \in A}$ in $B$. For $\alpha \in A$ we have $-\frac{1}{\lambda} x_\alpha \leq 0 \leq \frac{1}{f(x_\alpha)}x_\alpha \leq \frac{1}{\lambda}x_\alpha$, i.e.\ $\pm y_\alpha \leq \frac{1}{\lambda} x_\alpha$. This proves $y_\alpha \xrightarrow{o}0$. Thus since $B$ is order closed we obtain the contradiction $0\in B$ from Remark \ref{rem:basicorderconvergencestuff}(a). 
	\end{proof}
	
	The following proposition gives a necessary condition on $X$ to have an order closed base.
		
\begin{proposition} \label{pro:existslineartopologystBcompact_thenBorderclosed}
	Let $B$ be a base of $X_+$. If there exists a linear topology $\tau$ on $X$ such that $B$ is compact with respect to $\tau$, then $B$ is order closed.  	
\end{proposition}

\begin{proof}
	Let $f\colon X\to \mathbb{R}$ be the strictly positive linear functional on $X$ such that $B=\{x\in X_+ ;\, f(x)=1\}$. We use Remark \ref{rem:characterisationordercontinuityforpositivemaps}(b) to show that $f$ is order continuous. Let $(x_\alpha)_{\alpha \in A}$ be a map in $X$ such that $x_\alpha \downarrow 0$. Since $f$ is positive we obtain $(f(x_\alpha))_{\alpha \in A}$ to be a net in $X_+$ such that $f(x_\alpha)\downarrow$. To obtain a contradiction assume that $\lambda:=\inf\{f(x_\alpha);\, \alpha \in A\}>0$. For $\alpha \in A$ let $y_\alpha:=\frac{1}{f(x_\alpha)}x_\alpha$ to obtain a net $(y_\alpha)_{\alpha \in A}$ in $B$. Since $B$ is compact with respect to $\tau$ there is a subnet $(z_\gamma)_{\gamma \in C}$ of $(y_\alpha)_{\alpha \in A}$ such that $z_\gamma \xrightarrow{\tau}z\in B$. Thus there is a monotone map $h\colon C \to A$ such that $h[C]$ is majorizing in $A$ and such that $z_\gamma=y_{h(\gamma)}$ for every $\gamma \in C$. Since $(f(x_{h(\gamma)}))_{\gamma \in C}$ is a subnet of $(f(x_\alpha))_{\alpha \in A}$ we conclude $f(x_{h(\gamma)})\downarrow \lambda$. Thus Remark \ref{rem:basicorderconvergencestuff}(b) implies $f(x_{h(\gamma)})\xrightarrow{\tau_o(\mathbb{R})} \lambda$. Since $\tau_o(\mathbb{R})$ is the standard topology and $\tau$ is a linear topology on $X$ we conclude $x_{h(\gamma)}=f(x_{h(\gamma)})y_{h(\gamma)}=f(x_{h(\gamma)})z_\gamma \xrightarrow{\tau} \lambda z$. Furthermore we have $x_{h(\gamma)}\downarrow$. Thus Remark \ref{rem:basicsonordertopology}(b) shows $x_{h(\gamma)}\downarrow \lambda z$. We therefore obtain $\lambda z=0$ from $x_{h(\gamma)}\downarrow 0$. Since $z\in B$ this implies $\lambda=0$, a contradiction. 
\end{proof}

We next show that upper bounds of order closed bases of the cone are net catching. 

	\begin{lemma} \label{lem:descendingnetsfallunderorderclosedbases}
		Let $(x_\alpha)_{\alpha \in A}$ be a net in $X$ with $x_\alpha \downarrow 0$. 
		If $B$ is an order closed base of $X_+$
		then there is $b \in B$ and $\alpha \in A$ such that $x_\alpha \leq b$. 
	\end{lemma}
	
	\begin{proof}
	There is a strictly positive functional $f:X \rightarrow \R$ such that $B=\{x \in X_+;\, f(x)=1\}$. By Proposition \ref{pro:Borderclosedif_f_ordercontinuous} we obtain $f$ to be order continuous. Hence $f(x_\alpha)\downarrow 0$ by Remark \ref{rem:characterisationordercontinuityforpositivemaps}(b). Let $\alpha \in A$ be such that $f(x_\alpha)\leq 1$. If $f(x_\alpha)=0$, then the strict positivity of $f$ implies $x_\alpha=0$ and every $b\in B$ satisfies $x_\alpha \leq b$. If $f(x_\alpha)>0$, then $b := \frac{1}{f(x_\alpha)}x_\alpha$ satisfies $b \in B$ and $x_\alpha=f(x_\alpha) \frac{1}{f(x_\alpha)}x_\alpha \leq b$.
	\end{proof}

	\begin{proposition} \label{pro:Bleqxandorderclosedimpliesx_nc}
	If $B$ is an order closed base of $X_+$, then all upper bounds of $B$ are net catching.  
	\end{proposition}
	
	\begin{proof} Let $x\in X$ be an upper bound of $B$ and $(x_\alpha)_{\alpha \in A}$ a net in $X$ such that $x_\alpha \downarrow 0$. By Lemma \ref{lem:descendingnetsfallunderorderclosedbases} there are $\alpha \in A$ and $z \in B$ such that $x_\alpha \leq b$. Hence $x_\alpha \leq b \leq x$. 
	\end{proof}
	
	Combining Proposition \ref{pro:existslineartopologystBcompact_thenBorderclosed} with Proposition \ref{pro:Bleqxandorderclosedimpliesx_nc} we obtain the following Theorem.
	
	\begin{theorem} \label{the:necessarryconditionsfornetcatchingelements}
		Let $B$ be a base of $X_+$. If there exists a linear topology $\tau$ on $X$ such that $B$ is compact with respect to $\tau$, then all upper bounds of $B$ are net catching elements, i.e.\ interior points of $X_+$ with respect to the order topology. 
	\end{theorem}

\section{Examples of partially ordered vector spaces with net catching elements}
	In this Section we will apply Theorem \ref{the:necessarryconditionsfornetcatchingelements} to obtain that Archimedean finite dimensional partially ordered vector spaces and reflexive Banach spaces equipped with ice cream cones have net catching elements. 

\begin{proposition}\label{pro:finitedimensionalandArchimedeanspaceshavenetcatchinglements}
		Every nontrivial\footnote{A partially ordered vector space $X$ is said to be \emph{trivial}, if $X_+=\{0\}$.} Archimedean and finite dimensional partially ordered vector space has net catching elements.
	\end{proposition}
	
	\begin{proof} Note that a nontrivial partially ordered vector space $X$ has net catching elements if and only if its directed part $X_+-X_+$, equipped with the induced order from $X$, has net catching elements. Thus we can assume without lost of generality that $X$ is directed. Since $X$ is Archimedean we obtain\footnote{To obtain this combine \cite[Theorem 1.54]{CAD} with \cite[Corollary 3.8]{CAD}.} the existence of a base $B$ of $X_+$ such that $B$ is compact with respect to the standard topology. Since\footnote{To obtain this combine \cite[Lemma 2.5]{CAD} with \cite[Lemma 3.2]{CAD}.} $X$ is directed and finite dimensional there exists an order unit $u$ in $X$. We will now present, that there is $\lambda \in \mathbb{R}_{\geq 0}\setminus \{0\}$ such that $B\leq \lambda u$. Thus Theorem \ref{the:necessarryconditionsfornetcatchingelements} can be applied to show $\lambda u$ to be net catching. 
	
	Let $||\cdot||$ be a norm on $X$ such that the closed unit ball $\overline{B}_1(0)$ is the convex hull of a finite set $F\subseteq X$. Since $u$ is an order unit there is $\mu \in \mathbb{R}_{\geq 0}\setminus \{0\}$ such that $F\leq \mu u$, hence $\overline{B}_1(0)\leq \mu u$. Since $B$ is compact with respect to $||\cdot||$ there is $\kappa \in \mathbb{R}_{\geq 0}\setminus \{0\}$ such that $B\subseteq \overline{B}_\kappa(0)$. Thus letting $\lambda:=\kappa \mu$ we obtain $B\subseteq \{\kappa x;\, x \in\overline{B}_1(0)\} \subseteq \{\kappa x;\, x \leq \mu u \}\leq \lambda u$.
	\end{proof}

\begin{remark}
	(i) Note that the proof of the previous proposition shows that any order unit of an Archimedean finite dimensional partially ordered vector space is net catching. Indeed, we have shown that if $u$ is an order unit, then there is $\lambda \in \mathbb{R}_{\geq 0}\setminus \{0\}$ such that $\lambda u$ is net catching. Thus $u$ is net catching by Remark \ref{rem:ncelements}(iiib). 
	
	(ii) If $X$ is a directed Archimedean and finite dimensional partially ordered vector space, then the order units of $X$ are exactly the net catching elements. Indeed, any order unit is net catching by (i) and the reverse implication is shown in Proposition \ref{pro:nce_are_orderunits_XArchm_directed}.
	
	(iii) It remains an open question, whether every nontrivial finite dimensional partially ordered vector space has net catching elements.
\end{remark}

\begin{definition}
	Let $(X,||\cdot||)$ be a normed vector space. Let $f$ be a linear functional on $X$ such that $||f||:=\sup\{\frac{|f(x)|}{||x||};\, x \in X\}=1$ and $\epsilon \in (0,1)$. Then the \textbf{ice cream cone} with parameters $f$ and $\epsilon$ is the cone 
		\begin{align*}
			K_{f,\epsilon}:=\{x \in X;\,f(x) \geq \epsilon ||x||\}.
		\end{align*}
		The linear functional $f$ is called the \textbf{axis} and $\omega_\epsilon:=\operatorname{arccos}(\epsilon)$ is known as the \textbf{angle} of the ice cream cone $K_{f,\epsilon}$. 
\end{definition}

We will now show that any ice cream cone has net catching elements. 
	
	\begin{proposition} \label{pro:basicpropertiesoftheicecreamcone}
		Let $(X,||\cdot||)$ be a normed vector space and an $K_{f,\epsilon}$ an ice cream cone in $X$.
		\begin{itemize}
			\item[(i)] The cone $K_{f,\epsilon}$ is closed with a non-empty interior with respect to the norm topology.	Moreover this interior is given by $\{x \in X;\, f(x) > \epsilon ||x||\}$ and $K_{f,\epsilon}$ is generating.
			\item[(ii)] The boundary of $K_{f,\epsilon}$ is the set $\{x \in X;\,f(x)= \epsilon ||x||\}$. 
			\item[(iii)] The parameters of an ice cream cone are uniquely determined, i.e.\ two ice cream cones $K_{f,\epsilon}$ and $K_{g,\delta}$ satisfy $K_{f,\epsilon}=K_{g,\delta}$ if and only if $f=g$ and $\epsilon=\delta$. 
			\item[(iv)] The linear functional $f$ is $K_{f,\epsilon}$-strictly positive and therefore the closed convex set 
			\begin{align*}
				B_{f,\epsilon}:=\{x \in K_{f,\epsilon};\,f(x)=1\}
			\end{align*}
			is a base of $K_{f,\epsilon}$.
			\item[(v)] The partially ordered vector space $(X,K_{f,\epsilon})$ is Archimedean. 
		\end{itemize}
	\end{proposition}

\begin{proof}
	For (i),(ii),(iii) and (iv) see \cite[Theorem 2.52.]{CAD}. (v) can be obtained from (i) and \cite[Lemma 2.3(2)]{CAD}.
\end{proof}

	To show the existence of net catching elements in $(X,K_{f,\epsilon})$ we will apply Theorem \ref{the:necessarryconditionsfornetcatchingelements} to the base $B_{f,\epsilon}$.
	
	\begin{lemma} \label{lem:Basisoficecreamconeisniceintersection}
	Let $K_{f,\epsilon}$ be an ice cream cone in a normed space $(X,||\cdot||)$. Then $B_{f,\epsilon}= \overline{B}_{\frac{1}{\epsilon}}(0)\cap f^{-1}(\{1\})$.
	\end{lemma}
	
\begin{proof}
	\begin{align*}
	B_{f,\epsilon}&= \{x\in X;\, f(x)\geq \epsilon ||x|| \text{ and } f(x)=1\}\\ &=\{x\in X;\, \frac{1}{\epsilon}\geq ||x|| \text{ and } f(x)=1\}\\ &=\overline{B}_{\frac{1}{\epsilon}}(0)\cap f^{-1}(\{1\})
	\end{align*}
\end{proof}

	\begin{lemma} \label{lem:icecreamconebaseisweaklycompact}
	If $(X,||\cdot||)$ is a reflexive Banach space, then $B_{f,\epsilon}$ is weakly compact. 
	\end{lemma}
	
\begin{proof}
	If $X$ is a reflexive Banach space we know $\overline{B}_{\frac{1}{\epsilon}}(0)$ to be weakly compact. 
Thus $B$ is weakly compact as the intersection of the weakly closed preimage $f^{-1}(\{1\})$ with $\overline{B}_{\frac{1}{\epsilon}}(0)$.
\end{proof}

	\begin{lemma} \label{lem:Baseoficecreamconecanbemajorized}
		For each norm topological interior point $x$ of $K_{f,\epsilon}$ 
		there is $\lambda >0$ such that $B_{f,\epsilon} \leq \lambda x$. 
	\end{lemma}

\begin{proof} If $x$ is a norm topological interior point $K_{f,\epsilon}$, then there is $\kappa >0$ such that $\overline{B}_{\kappa}(x)\subseteq K_{f,\epsilon}$. Define $\lambda:=\frac{1}{\epsilon \kappa}$ and let $b \in B_{f,\epsilon}$. From Lemma \ref{lem:Basisoficecreamconeisniceintersection} we know $||b||\leq \frac{1}{\epsilon}$. Thus
		 \begin{align*}
		 	||x-\epsilon \kappa (\lambda x-b)||
		 			= ||\epsilon \kappa b||
		 			= \epsilon \kappa ||b||
		 			\leq \kappa.
		 \end{align*}
		 Hence $\epsilon \kappa(\lambda x-b) \in \overline{B}_{\kappa}(x) \subseteq K_{f,\epsilon}$. Since $K_{f,\epsilon}$ is a cone
		 we obtain $\lambda x-b \in K_{f,\epsilon}$, i.e.\ $b \leq \lambda x$.
\end{proof}
	
\begin{proposition} \label{pro:icecramconehasmajorizingorderunitsifXisreflexivebanachspace}
		Let $X$ be a reflexive Banach space. Every norm topological interior point of $K_{f,\epsilon}$ is net catching. 
\end{proposition}

\begin{proof}
	From Theorem \ref{the:necessarryconditionsfornetcatchingelements} and Lemma \ref{lem:icecreamconebaseisweaklycompact} we obtain, that any upper bound of $B_{f,\epsilon}$ is a net catching element. Furthermore if $x$ is an norm topological interior point of $K_{f,\epsilon}$ we have shown in  Lemma \ref{lem:Baseoficecreamconecanbemajorized}, that there is $\lambda>0$ such that $\lambda x$ is an upper bound of $B_{f,\epsilon}$, hence net catching. Thus $x$ is net catching by Remark \ref{rem:ncelements}(iiib).
\end{proof}
	
\section{Net catching elements in Archimedean vector lattices}

In Proposition \ref{pro:finitedimensionalandArchimedeanspaceshavenetcatchinglements} we have seen that an Archimedean vector space has net catching elements, whenever it is finite dimensional. We will show in this section that the converse is valid in Archimedean vector lattices. We will use the following representation theorem for Dedekind complete vector lattices with order units. 

	\begin{proposition} \label{pro:representationCOmega}
		For every Dedekind complete vector lattice $X$ with an order unit $u$ there exists an extremely disconnected\footnote{A topological space $\Omega$ is called \emph{extremely disconnected}, if the closure of every open set is open.} and compact Hausdorff space $(\Omega,\tau)$ and a linear order isomorphism $J:X\to C(\Omega)$ such that $J(u)(\omega)=1$ for every $\omega \in \Omega$, where $C(\Omega)$ denotes the partially ordered vector space of all bounded and continuous maps $f\colon \Omega \to \mathbb{R}$.
	\end{proposition}
	
	\begin{proof}
		Let $\mathbb{R}_\infty:=\mathbb{R}\cup\{-\infty,\infty\}$, where $\infty$ and $-\infty$ are symbols not being used in $\mathbb{R}$. If $\tau_\mathbb{R}$ is the standard topology on $\mathbb{R}$, then $B_\infty:=\tau_\mathbb{R}\cup\{(\lambda,\infty)\cup \{\infty\};\, \lambda \in \mathbb{R}\}\cup\{(-\infty,\lambda)\cup \{-\infty\};\, \lambda \in \mathbb{R}\}$ is a base  of a topology $\tau_\infty$ on $R_\infty$. We equip $R_\infty$ with this topology. If $(\Omega,\tau)$ is an extremely disconnected compact Hausdorff space we define $C_\infty(\Omega)$ as the set of all continuous functions $f\colon \Omega \to \mathbb{R}_\infty$ such that $[\{-\infty,\infty\}]f$ is nowhere dense in $\Omega$. In \cite[Lemma V.2.1 ff.]{Vulikh67} it is shown that for $f,g\in C_\infty(\Omega)$ there is a unique $h\in C_\infty(\Omega)$ such that $h(\omega)=f(\omega)+g(\omega)$ for every $\omega \in [\mathbb{R}]f\cap [\mathbb{R}]g$. This defines an addition on $C_\infty(\Omega)$. For $\lambda \in \mathbb{R}$ we set $\lambda(\pm \infty):=\pm \infty$, if $\lambda>0$, $\lambda(\pm \infty):=0$, if $\lambda=0$ and $\lambda(\pm \infty):=\mp \infty$, if $\lambda<0$. Hence $(\lambda f)(\omega):=\lambda f(\omega)$ for every $\omega \in \Omega$ and $f \in C_\infty(\Omega)$ defines a scalar multiplication in $C_\infty(\Omega)$. If we furthermore order $C_\infty(\Omega)$ component wise we obtain a Dedekind complete vector lattice, as presented in \cite[Theorem V.2.2]{Vulikh67}.
		
In \cite[Theorem V.4.1]{Vulikh67} it is shown\footnote{Note that the Theorem in \cite[Theorem V.4.1]{Vulikh67} is stated for Dedekind complete vector lattices $X$ with 'units', i.e.\ elements $u\in X$ such that $x\wedge u>0$, whenever $x>0$. Nevertheless each order unit is a 'unit'. Furthermore in \cite[Theorem V.4.1]{Vulikh67} $Y$ is assumed to be a 'foundation', i.e.\ an ideal with further properties. Nevertheless from \cite[Definition II.9.1]{Vulikh67} we obtain that every ideal is a vector subspace. Furthermore note that a 'algebraically and lattice isomorphism' is nothing but a linear order isomorphism.} that for every Dedekind complete vector lattice $X$ with order unit $u$ there exists a extremely disconnected and compact Hausdorff space $(\Omega,\tau)$, a vector subspace $Y$ of $C_\infty(\Omega)$ and a linear order isomorphism $J\colon X\to Y$ such that $C(\Omega)\subseteq Y$ and $J(u)(\omega)=1$ for every $\omega \in \Omega$. To prove the statement it is sufficient to show $Y=C(\Omega)$. Let $f\in Y$ and note that $J(u)$ is an order unit in $Y$. Thus there is $\lambda \in \mathbb{R}$ such that $\pm f\leq \lambda J(u)$. For $\omega \in \Omega$ we therefore obtain $-\lambda =-\lambda J(u)(\omega)\leq f(\omega) \leq \lambda J(u)(\omega)=\lambda$, hence $f\in C(\Omega)$. This proves $Y=C(\Omega)$. 
	\end{proof}

	\begin{lemma} \label{lem:constructionofweirdnet}
		Let $P$ be a lattice, $Q$ a partially ordered set and $f \colon P \to Q$ an order embedding such that $f[P]$ is order dense in $Q$. Let $(y_\alpha)_{\alpha \in A}$ be a net in $Q$ such that $y_\alpha \downarrow f(x)$ for $x\in P$. If $B:=\{v\in P;\, \exists \alpha \in A: f(v)\geq y_\alpha\}$ is equipped with the reversed order from $P$, then $B$ is directed and $\inf B=x$. Thus $x_\beta:=\beta$ for all $\beta \in B$ defines a net in $P$ with $x_\beta \downarrow x$.
	\end{lemma}
	
\begin{proof}
	See \cite[Lemma 4.9]{Paper1}.
\end{proof}
	
	\begin{proposition} \label{pro:orderembeddingsandnetcatchingelements}
		Let $P$ be a lattice, $Q$ a partially ordered set, $f\colon P\to Q$ an order embedding such that $f[P]$ is order dense in $Q$, and $x,y\in P$. Then $y$ is downward net catching with respect to $x$ if and only if $f(y)$ is downward net catching with respect to $f(x)$. A similar statement is true about upward net catching elements.
		
	Let $X$ be a vector lattice, $Y$ a partially ordered vector space, $f\colon X\to Y$ a linear order embedding such that $f[X]$ is order dense in $Y$, and $x,y\in X$. Then $x$ is net catching if and only if $f(x)$ is net catching.
	\end{proposition}
	
\begin{proof} If $(x_\alpha)_{\alpha \in A}$ is a net in $P$ with $x_\alpha\downarrow x$, then $f(x_\alpha)\downarrow f(x)$ by Remark \ref{rem:characterisationordercontinuityforpositivemaps}(a). Thus if $f(y)$ is downward net catching with respect to $f(x)$ there is $\alpha \in A$ with $f(x_\alpha)\leq f(y)$, i.e.\ $x_\alpha \leq y$. To show the reverse implication of the first statement let $(y_\alpha)_{\alpha \in A}$ be a net in $Q$ such that $y_\alpha \downarrow f(y)$. Let $(x_\beta)_{\beta \in B}$ be the net defined in Lemma \ref{lem:constructionofweirdnet} and note that $x_\beta\downarrow x$. Since $y$ is assumed to be net catching with respect to $x$ there exists $\beta \in B$ such that $x_\beta \leq y$. Hence there is $\alpha \in A$ with $y_\alpha \leq f(\beta)=f(x_\beta) \leq f(y)$. This proves $f(y)$ to be downward net catching with respect to $f(x)$. The remaining statements are a direct consequence of the first one. 
\end{proof}

	\begin{proposition} \label{pro:Dedekindcomplete_Rieszspace_existence_netcatchingelements_implies_finitedimensional}
	Let $X$ be a  Dedekind complete vector lattice. If $X$ has net catching elements, then it is finite dimensional.
	\end{proposition}
	
	\begin{proof}
		Note that $X$ is directed and Archimedean. Thus if $X$ has no order units, then $X$ can not have net catching elements by Proposition \ref{pro:nce_are_orderunits_XArchm_directed}. We therefore assume without lost of generality that $X$ has an order unit $u$. By Proposition \ref{pro:ncexist_thenorderunitsare_nc} it is sufficient to show that $u$ is not net catching. Since $u$ is an order unit Proposition \ref{pro:representationCOmega} shows the existence of an extremely disconnected compact Hausdorff space $(\Omega,\tau)$ and a linear order isomorphism $J\colon X\to C(\Omega)$ such that $J(u)(\omega)=1$ for every $\omega \in \Omega$. We will now construct a net $(x_\alpha)_{\alpha \in A}$ in $C(\Omega)$ such that $x_\alpha \downarrow 0$, but $x_\alpha\leq J(u)$ for no $\alpha \in A$. Thus $J(u)$ is not net catching and Proposition \ref{pro:orderembeddingsandnetcatchingelements} therefore implies $u$ to be not net catching. 
		
Since $C(\Omega)$ is infinite dimensional we know $\Omega$ to be infinite. Thus the compactness of $\Omega$ implies the existence of a not isolated\footnote{A point $\omega$ in a topological space $\Omega$ is called \emph{isolated}, if $\{\omega\}$ is open and closed.} point $\delta \in \Omega$. Let $A:=\{M \subseteq \Omega;\,\delta \in M \text{ and } M \text{ open and closed}\}$ be ordered by reversed inclusion. Since the intersection of two sets in $A$ is contained in $A$ we obtain that $A$ is directed. For $\alpha \in A$ let $M_\alpha:= \alpha$ and $\frac{1}{2}x_\alpha$ be the characteristic function of $M_\alpha$. Since $(M_\alpha)_{\alpha \in A}$ is a decreasing net of open and closed sets with respect to set inclusion we obtain $(x_\alpha)_{\alpha \in A}$ to be a decreasing net in $C(\Omega)$. Note furthermore that $x_\alpha(\delta)=2>1=J(u)(\delta)$ for every $\alpha \in A$, hence $x_\alpha \leq J(u)$ for no $\alpha \in A$. It is left to show that $\inf\{x_\alpha;\, \alpha \in A\}=0$. Clearly $0$ is a lower bound of $\{x_\alpha;\, \alpha \in A\}$. Let $z\in C(\Omega)$ be another lower bound of this set and $\omega \in \Omega\setminus \{\delta\}$. Since $\Omega$ is Hausdorff  there are disjoint $U,O\in \tau$ such that $\delta \in U$, $\omega
 \in O$. Let $M$ be the closure of $U$ and note that $M\in \tau$, since $\Omega$ is extremely disconnected. Furthemore $\delta \in M$ and $\omega \notin M$. Hence $M \in A$ and we obtain $z(\omega)\leq x_M(\omega)=0$. We have shown $z(\omega)\leq 0$ for every $\omega \in \Omega \setminus\{\delta\}$. To obtain a contradiction assume $z(\delta)>0$. Hence $\{\delta\}=[(0,\infty) \cup \{\infty\}]z$ is open and closed, a contradiction to the assumption that $\delta$ is isolated. This proves $z\leq 0$ and we have shown $x_\alpha \downarrow 0$. 
	\end{proof}

As a combination of Proposition \ref{pro:finitedimensionalandArchimedeanspaceshavenetcatchinglements}, Proposition \ref{pro:orderembeddingsandnetcatchingelements} and Proposition \ref{pro:Dedekindcomplete_Rieszspace_existence_netcatchingelements_implies_finitedimensional}
 we obtain the following. 

\begin{theorem} \label{the:Characterization_existencence_for_Archimedean_VL}
	Let $X$ be an Archimedean vector lattice. $X$ has net catching elements if and only if $X$ is finite dimensional. 
\end{theorem}

\begin{proof} In Proposition \ref{pro:finitedimensionalandArchimedeanspaceshavenetcatchinglements}
	we have seen that finite dimensional Archimedean partially ordered vector spaces have net catching elements. To show the converse assume $X$ to be infinite dimensional and let $(X^\delta,J)$ be a Dedekind completion of $X$. Clearly $X^\delta$ is infinite dimensional. Thus Proposition \ref{pro:Dedekindcomplete_Rieszspace_existence_netcatchingelements_implies_finitedimensional}
shows that $X^\delta$ has no net catching elements. Hence $X$ can not have net catching elements by Proposition \ref{pro:orderembeddingsandnetcatchingelements}.
\end{proof}

\begin{remark}
Note that Example \ref{exa:lexicographicsequencespacenetcatchingelements} provides an infinite dimensional vector lattice with net catching elements. Thus the Archimedean property can not be omitted in Theorem \ref{the:Characterization_existencence_for_Archimedean_VL}.

Furthermore, as seen in Proposition \ref{pro:icecramconehasmajorizingorderunitsifXisreflexivebanachspace}, every ice cream cone in a reflexive Banach space has net catching elements. Thus there are infinite dimensional Archimedean partially ordered vector spaces with net catching elements. The assumption on $X$ to be a vector lattice can therefore not be omitted in Theorem \ref{the:Characterization_existencence_for_Archimedean_VL}. 
\end{remark} 

\section{The Minkowski norm and the order bound topology}

In this section let $X$ denote a partially ordered vector space. In \cite{Namioka57} and \cite{CAD} the concept of the order bound topology is introduced and studied. This concept differs from the concept of order topology as we will see. We will show that the order units play an analogous role for the order bound topology as the net catching elements for the order topology. We will furthermore provide cases in which a relation between these topologies can be obtained and show that net catching sets for $0$ contain circled sets, and absorb all order intervals if $X$ is Archimedean and directed. 

\begin{definition}
	Let $\mathcal{B}_{ob}(X)$ (or simply $\mathcal{B}_{ob}$) be the set of all subsets $V \subseteq X$ such that $V$ is convex, circled\footnote{A subset $M\subseteq V$ of a vector space $V$ is called \emph{circled}, if $\lambda M\subseteq M$ for every $\lambda \in [-1,1]$. If $M,N\subseteq V$, then $M$ is said to \emph{absorb} $N$, if there is $\mu >0$ such that $\lambda N\subseteq M$ for every $\lambda \in [-\mu,\mu]$.} and absorbs all order intervals. Note that $\mathcal{B}_{ob}$ is a base at zero for a linear topology $\tau_{ob}(X)$ (or simply $\tau_{ob}$). This topology is called the \emph{order bound topology}.
\end{definition}

\begin{remark}
	There exists a non empty order bounded and $\tau_{ob}$-open set in $X$ if and only if $\tau_{ob}\operatorname{-int}(X_+)\neq \emptyset$. This can be proven similar to Proposition \ref{pro:characterisation_existenceorderopenandorderbddsets}.
\end{remark}
		
\begin{definition}
	Assume $X$ to contain an order unit $u$. For $x\in X$ define the \emph{$u$-norm} as $||x||_u:=\inf\{\lambda > 0;\,\pm x \leq \lambda u\}$.
\end{definition}
	
\begin{remark}
The order bound topology $\tau_{ob}(X)$ is introduced and studied in \cite{Namioka57}. Further research on this subject can be found in \cite[2.8]{CAD}, where the order bound topology is referred to as 'order topology'. The definition of the $u$-norm can be found in \cite[2.7]{CAD}.
\end{remark}

\begin{remark} \label{rem:on_unorm}
Let $u$ and $u'$ be order units in $X$. 
\begin{itemize}
	\item[(i)] The map $||\cdot||_u$ is 
	a norm on $X$. Furthermore the corresponding closed unit Ball satisfies $\overline{B}_1(0)=[-u,u]$ and $X_+$ is norm closed with respect to $||\cdot||_u$. This is shown in \cite[Theorem 2.55.]{CAD}. 
	\item[(ii)] The topology generated by $||\cdot||_u$ and the order bound topology agree. This result can be found in \cite[Theorem 2.71]{CAD}. In particular the norms $||\cdot||_u$ and  $||\cdot||_{u'}$ are equivalent.
\end{itemize}
\end{remark}

The following proposition shows that the elements of $\mathcal{B}_{ob}$ and order units play a similar role as net catching sets for $0$ and net catching elements repectively in Proposition \ref{pro:nc_elements_and_nc_sets}
		
\begin{proposition} \label{pro:orderunitsgiveneighbourhoodsof0}
	Let $u\in X$. Then $u$ is an order unit if and only if $[-u,u] \in \mathcal{B}_{ob}$.
\end{proposition}

\begin{proof}
First assume $u$ to be an order unit. Clearly $[-u,u]$ is convex and circled. For $x,y\in X$ with $x\leq y$ there are $m,n\in \mathbb{N}$ such that $-x \leq m u$ and $y \leq n u$. Let $\lambda >0$ such that $m \leq \frac{1}{\lambda}$ and $n \leq \frac{1}{\lambda}$ to obtain $-u\leq \lambda x\leq \lambda y\leq u$. Since $[-u,u]$ is circled this proves $[-u,u]$ to absorb $[x,y]$. 

To show the converse implication assume $[-u,u]\in \mathcal{B}_{ob}$ and let $x\in X$. Since $[-u,u]$ absorbs all order intervals there is $\lambda >0$ such that $\lambda[x,x]\subseteq [-u,u]$. Let $n \in \mathbb{N}_{\geq \frac{1}{\lambda}}$ to obtain $x \leq n u$.
\end{proof} 
		
In Theorem \ref{the:interiorofconeare_nc_elements} we have seen that the interior of the cone with respect to $\tau_o$ can be characterized as the set of net catching elements.  The following analogue is valid. 

	\begin{proposition} \label{pro:orderboundtopologicalinteriorofconeareorderunits}
	The order bound topological interior of the cone is the	set of all order units, i.e.\ $\tau_{ob}\operatorname{-int}(X_+)=\{u \in X;\,u \text{ is an order unit}\}$.
\end{proposition}
\begin{proof} Let $u \in \tau_{ob}\text{-int}(X_+)$ and let $x\in X$. There is $V \in \mathcal{B}_{ob}(X)$ such that $V+u \subseteq X_+$. Since $V$ absorbs order intervals, there is $\mu>0$ such that $\mu [-x,x] \subseteq V$, hence $-\mu x+u \in \mu [-x,x]+u \subseteq V+u \subseteq X_+$. Let $n\in \mathbb{N}_{\geq \frac{1}{\mu}}$ to obtain $x\leq n u$. Using Proposition \ref{pro:orderunitsgiveneighbourhoodsof0} a similar argument as in the proof of Theorem \ref{the:interiorofconeare_nc_elements} gives $\{u \in X;\,u \text{ is an order unit}\}\in \tau_{ob}$.
\end{proof}

\begin{lemma} \label{lem:infimumoverorderunits}
	Assume $X$ to contain order units and to be Archimedean. Let $U$ be the set of all order units of $X$. Then $U$ is downward directed and $\inf U=0$.
\end{lemma}

\begin{proof}
Let $u,v\in U$. There is $n \in \mathbb{N}$ such that $v \leq n u$. Thus $\frac{1}{n}u \leq v$ and $\frac{1}{n}u \leq u$. Since $\frac{1}{n}u\in U$ we obtain $U$ to be downward directed. 

Clearly $0$ is a lower bound of $U$. Let $z$ be another lower bound of $U$ and $u \in U$. Since $X$ is Archimedean we obtain $\frac{1}{n}u\downarrow 0$. Since $\left(\frac{1}{n}u\right)_{n \in \N}$ is a net in $U$ we conclude $z \leq 0$. This shows $\inf U=0$.
\end{proof}
			
\begin{proposition} \label{pro:orderunitsandArchimedeanimplynormconvergenceimplieso3convergence}
	Let $X$ be Archimedean, $u\in X$ an order unit and $(x_\alpha)_{\alpha \in A}$ a net in $X$. Whenever $x_\alpha \xrightarrow{||\cdot||_u} x\in X$, then also $x_\alpha \xrightarrow{o} x$.
\end{proposition}

\begin{proof}
	Without lost of generality assume that $x=0$. Let $B$ be the set of all order units, equipped with the reversed order from $X$. Define $\check{x}_\beta:=\beta$ for $\beta \in B$ to obtain a net $(\check{x}_\beta)_{\beta \in B}$. Lemma \ref{lem:infimumoverorderunits} implies that $\check{x}_\beta \downarrow 0$. To construct a map $\eta: B \rightarrow A$ let $\beta \in B$. Since $\beta$ is an order unit of $X$ Remark \ref{rem:on_unorm}(ii) implies the equivalence of the norms $||\cdot||_u$ and $||\cdot||_{\beta}$ and we are able to conclude $x_\alpha \xrightarrow{||\cdot||_\beta}0$. Thus there is $\alpha_\beta \in A$ such that $\inf \{\lambda >0;\,\pm x_\alpha \leq \lambda \beta\}=||x_\alpha||_\beta\leq \frac{1}{2}$ for $\alpha \in A_{\geq \alpha_\beta}$. We obtain $\pm x_\alpha \leq \beta=\check{x}_\beta$ for every $\alpha \in A_{\geq \alpha_\beta}$. Hence $\eta(\beta):=\alpha_\beta$ allows to conclude $x_\alpha \xrightarrow{o} 0$.
\end{proof} 

\begin{remark} 
	Note that the net $(\check{x}_\beta)_{\beta \in B}$ is independent from $u$ in the construction in the proof of Proposition \ref{pro:orderunitsandArchimedeanimplynormconvergenceimplieso3convergence}.
\end{remark}

The following theorem gives a relationship between the order topology and the order bound topology.

\begin{theorem} \label{the:To_and_Tob}
	If $X$ is a Archimedean and if there are order units in $X$, then the order bound topology is finer than the order topology, i.e.\ $\tau_o\subseteq \tau_{ob}$.
	If there are net catching elements in $X$, then the order topology is finer than the order bound topology, i.e.\ $\tau_{ob} \subseteq  \tau_o$.	
\end{theorem}

\begin{proof}
Let $u$ be an order unit in $X$. If $(x_\alpha)_{\alpha \in A}$ is a net in $X$ with $x_\alpha \xrightarrow{\tau_{ob}}x\in X$, then $x_\alpha \xrightarrow{||\cdot||_u} x\in X$ by Remark \ref{rem:on_unorm}(ii). Furthermore Proposition \ref{pro:orderunitsandArchimedeanimplynormconvergenceimplieso3convergence} implies $x_\alpha \xrightarrow{o} x$. Thus $x_\alpha \xrightarrow{\tau_o} x$ by Remark \ref{rem:basicorderconvergencestuff}(b). This proves the first statement. 
	
	To show the second statement we will prove that every $U\in \mathcal{B}_{ob}$ is net catching for $0$. Let $u$ be a net catching element. Since $U$ absorbs all order intervals there is $\lambda >0$ such that $[-\lambda u,\lambda u]\subseteq U$. By Remark \ref{rem:ncelements}(iiib) $\lambda u$ is net catching as well. Thus Proposition \ref{pro:nc_elements_and_nc_sets} shows $[-\lambda u,\lambda u]$ to be net catching for $0$. From $[-\lambda u,\lambda u]\subseteq U$ we therefore obtain $U$ to be net catching for $0$. 
\end{proof}

\begin{remark} 
	(i) Note that the first inclusion in Theorem \ref{the:To_and_Tob} can be strict and that the existence of net catching elements can not be omitted in the second statement. Indeed, in Example \ref{exa:l_inftyhasnonetcatchingelements} we have seen, that $l_\infty$ with the component wise order has order units, but no net catching element. Thus Proposition \ref{pro:orderboundtopologicalinteriorofconeareorderunits} implies the interior of the cone with repsect to $\tau_{ob}$ to be nonempty. Nevertheless the interior of the cone with respect to $\tau_o$ is empty by Theorem \ref{the:interiorofconeare_nc_elements}, which implies $\tau_o\neq \tau_{ob}$. Since $\tau_o\subseteq \tau_{ob}$ by Theorem \ref{the:To_and_Tob} we obtain $\tau_{ob}\not \subseteq \tau_o$. 
	
	(ii) The second inclusion in Theorem \ref{the:To_and_Tob} can be strict as well and the Archimedean property can not be omitted in the first statement. Indeed, in Example \ref{exa:lexicographicsequencespacenetcatchingelements} we have seen a partially ordered vector space $X$ with net catching elements and that $X$ contains order units which are not net catching. Thus $\tau_o\operatorname{-int}(X_+)\neq \tau_{ob}\operatorname{-int}(X_+)$ by Theorem \ref{the:interiorofconeare_nc_elements} and Proposition \ref{pro:orderboundtopologicalinteriorofconeareorderunits}. Hence $\tau_o\not \subseteq \tau_{ob}$ by Theorem \ref{the:To_and_Tob}.
	
	(iii) It is an open question, whether the assumption on the existence order units can be omitted in the first statement in Theorem \ref{the:To_and_Tob}.
\end{remark}

In the remainder of this section we show for Archimedean and directed partially ordered vector spaces 
that the set of all net catching sets for $0$ that are circled and absorb all order intervals behave similarly as the elements of $\mathcal{B}_{ob}$. 
	In particular we establish that every net catching set for $0$ contains a circled net catching set for $0$ that absorbs all order intervals. 

\begin{lemma}\label{lem:circlednetcatchingsets}
	Let $U\subseteq X$ be a net catching set for $0$. Define $I(U):=\{x\in U;\, [-x,x]\subseteq U\}$ and $V(U):=\bigcup_{x\in I(U)}[-x,x]$. Then $V(U)$ is a circled 
	net catching set for $0$.
\end{lemma}

\begin{proof} 
	$V(U)$ is circled as the union of circled sets. 
	To show that $V(U)$ is net catching for $0$  we will use the characterization of net catching sets for $0$ in partially ordered vector spaces, which was presented in Remark \ref{rem:basicsonordertopology}(c). Let $(x_\alpha)_{\alpha\in A}$ be a net in $X$ with $x_\alpha \downarrow 0$. As $U$ is a net catching set for $0$ there is $\alpha \in A$ with $[-x_\alpha,x_\alpha]\subseteq U$. Thus $x_\alpha\in I(U)$ and we obtain $[-x_\alpha,x_\alpha]\in V(U)$.	
\end{proof}

\begin{remark}
	The construction of the set $V(U)$ is given in \cite{Imh} and it is shown that $V(U)$ is a net catching set for $0$, whenever $U$ is a net catching set for $0$. In this thesis this construction is used to show that every net catching set for $0$ contains a set $V$ that is a net catching set for $0$ with the property that for every $x\in V\cap X_+$ there is $[-x,x]\subseteq V$. 
\end{remark}

\begin{lemma} \label{lem:circledandnetcatchingimpliesintervalabsorbtion}
	If $X$ is Archimedean and directed and $V$ is a circled net catching set for $0$, then $V$ absorbs\footnote{A subset $M$ of $X$ is said to \emph{absorb} $N\subseteq X$, if there exists $\lambda_0>0$ such that $\lambda n\in M$ for all $\lambda \in [-\lambda_0,\lambda_0]$ and $n\in N$.} all order intervals. 
\end{lemma}

\begin{proof} 
	Let $x,y\in X$ with $x\leq y$. As $V$ is circled it is sufficient to show the existence of $\lambda_0$ such that $\lambda_0 [x,y]\subseteq V$. Since $X$ is directed there is $w\in X$ with $w\geq -x$ and $w\geq y$. As $X$ is Archimedean we have $\frac{1}{n}w\downarrow 0$ and we obtain the existence of $n\in \mathbb{N}$ with $\frac{1}{n}[x,y]\subseteq [-\frac{1}{n} w,\frac{1}{n}w]\subseteq V$ as $V$ is a net catching set for $0$. 
\end{proof}

\begin{proposition} \label{pro:Xnetcatchingsetscircledabsorbtion}
	Every net catching set for $0$ contains a circled net catching set for $0$. If $X$ is Archimedean and directed, then every net catching set for $0$ absorbs all order intervals. 
\end{proposition}

\begin{proof}
	As $V(U)\subseteq U$ for every $U\subseteq X$ the first statement follows from Lemma \ref{lem:circlednetcatchingsets}. If a subset of a set $U$ absorbs all order intervals, then so does $U$. Thus the second statement follows from the first statement and Lemma \ref{lem:circledandnetcatchingimpliesintervalabsorbtion}.
\end{proof}

\section{Order topology in partially ordered vector spaces with net catching elements}
In this section let $X$ be a partially ordered vector space. We will see that under the assumption of the existence of net catching elements in $X$ the order topology $\tau_o$ is normable if and only if $X$ is Archimedean and directed. Furthermore we will obtain under this assumption that the order topological convergence is equivalent to $o$-convergence.

\begin{lemma} \label{lem:Xarchimedeananddirectedif1overnxordertopologconvto0}
	$X$ is Archimedean and directed if and only if $\left(\frac{1}{n}\right)_{n \in \mathbb{N}}$ satsifies $\frac{1}{n} x \xrightarrow{\tau_o} 0$ for every $x\in X$. 
\end{lemma}

\begin{proof}
	This statement is proven in \cite[Lemma 8.1]{Paper1}.
\end{proof}

\begin{theorem} \label{the:ArchimedeananddirectedandnetcatchingthenTobisTo}
	If $X$ contains net catching elements, then the following statements are equivalent.
			\begin{itemize}
				\item[(i)] $\tau_o$ is a linear topology.
				\item[(ii)] $\tau_o=\tau_{ob}$.
				\item[(iii)] $X$ is Archimedean and directed.				
			\end{itemize}
In particular if one of these statements is valid, then each net catching element $x\in X$ is an order unit of $X$ and the norm $||\cdot||_x$ generates $\tau_o$.
\end{theorem}

\begin{proof} Since $\tau_{ob}$ is a linear topology clearly (ii) implies (i). If $\tau_o$ is a linear topology, then $\left(\frac{1}{n}\right)_{n \in \mathbb{N}}$ satsifies $\frac{1}{n} x \xrightarrow{\tau_o} 0$ for every $x\in X$. Thus we obtain $X$ to be Archimedean and directed from Lemma \ref{lem:Xarchimedeananddirectedif1overnxordertopologconvto0} and (i) implies (iii). If $X$ is Archimedean and directed, then every net catching element is an order unit by Proposition \ref{pro:nce_are_orderunits_XArchm_directed}. Thus Theorem \ref{the:To_and_Tob} yields $\tau_o=\tau_{ob}$ and we have shown that (ii) follows from (iii). 
\end{proof}

\begin{remark}
	Let $X$ be finite dimensional. Then $\tau_o$ is the standard topology\footnote{The \emph{standard topology} on a finite dimensional vector space $X$ is the only linear topology on $X$.} if and only if $X$ is Archimedean and directed. Indeed, if $\tau_o$ is the standard topology, i.e.\ linear, then we obtain $X$ to be Archimedean and directed from Lemma \ref{lem:Xarchimedeananddirectedif1overnxordertopologconvto0}. If $X$ is assumed to be Archimedean and directed, then $X$ has net catching elements by Proposition \ref{pro:finitedimensionalandArchimedeanspaceshavenetcatchinglements}. Thus Theorem \ref{the:ArchimedeananddirectedandnetcatchingthenTobisTo} implies $\tau_o$ to be linear, i.e.\ the standard topology.  
\end{remark}

In partially ordered vector spaces with net catching elements we can characterise the order topological convergence as $o$-convergence. 

\begin{corollary} \label{cor:ordertopconvergenceando3}
	If $X$ is a Archimedean and directed and contains net catching elements, then for every net $(x_\alpha)_{\alpha \in A}$ in $X$ and every $x\in X$ we have $x_\alpha \xrightarrow{o}x$ if and only if $x_\alpha \xrightarrow{\tau_o} x$. 
\end{corollary}

\begin{proof}
	If $x_\alpha \xrightarrow{o}x$, then also $x_\alpha \xrightarrow{\tau_o} x$ by Remark \ref{rem:basicorderconvergencestuff}(b). To show the converse implication let $x_\alpha \xrightarrow{\tau_o} x$ and $u$ a net catching element of $X$. Proposition \ref{pro:nce_are_orderunits_XArchm_directed} shows that $u$ is an order unit of $X$. By Theorem \ref{the:ArchimedeananddirectedandnetcatchingthenTobisTo} we therefore obtain $x_\alpha \xrightarrow{||\cdot||_u} x$. Hence $x_\alpha \xrightarrow{o}x$ by Proposition \ref{pro:orderunitsandArchimedeanimplynormconvergenceimplieso3convergence}. 
\end{proof}

\begin{remark}\label{rem:ImhoffsFehler}
	In \cite{Imh} the equivalence of the $o$-convergence and of 	the order topological convergence is proven with a different technique and stated in \cite[Corollary 4.11]{Imh} for Archimedean partially ordered vector spaces with order units. 
	Nevertheless in the proof the set of all order open, order bounded and $0$ containing sets is used as an index set for a net. In Proposition \ref{pro:characterisation_existenceorderopenandorderbddsets} and Theorem \ref{the:interiorofconeare_nc_elements} it is shown, that this index set is non empty if and only if $X$ has net catching elements. Therefore the content of \cite[Corollary 4.11]{Imh} is the same as Corollary \ref{cor:ordertopconvergenceando3} and one should be aware of this implicit assumption in the result \cite[Corollary 4.11]{Imh}.
\end{remark}

\section{The order bound topology of an ice cream cone}

In Proposition \ref{pro:icecramconehasmajorizingorderunitsifXisreflexivebanachspace} we have seen that in reflexive Banach spaces $(X,||\cdot||)$ each ice cream cone $K_{f,\epsilon}$ contains net catching order units. Note that $(X,K_{f,\epsilon})$ is Archimedean and directed by Proposition \ref{pro:basicpropertiesoftheicecreamcone}. Hence $\tau_{ob}=\tau_o$ by Theorem \ref{the:ArchimedeananddirectedandnetcatchingthenTobisTo} and we know that $||\cdot||_x$ induces $\tau_o$ for every net catching order unit $x\in X$. We will show that $||\cdot||$ and $||\cdot||_x$ are equivalent and thus $\tau_o$ is the norm topology of $X$.
		
\begin{proposition}
	Let $(X,||\cdot||)$ be a reflexive Banach space and equip $X$ with an ice cream cone $K_{f,\epsilon}$. Then the topology induced by $||\cdot||$ is the order topology $\tau_o$ generated by $K_{f,\epsilon}$, and agrees with the order bound topology $\tau_{ob}$.
\end{proposition}

\begin{proof} Let $x\in X$ be a net catching order unit. It is left to show that $||\cdot||_x$ and $||\cdot||$ are equivalent. In Remark \ref{pro:icecramconehasmajorizingorderunitsifXisreflexivebanachspace} we have seen, that every interior point $y$ of $K_{f,\epsilon}$ with respect to $||\cdot||$ is net catching. Since $||\cdot||_x$ and $||\cdot||_y$ are equivalent we can assume without lost of generality that $x\in ||\cdot||\operatorname{-int}(K_{f,\epsilon})$. Thus there is $\kappa >0$ such that $\overline{B}_\kappa^{||\cdot||}(x)\subseteq K_{f,\epsilon}$. Furthermore let $\lambda:=\frac{2}{\epsilon}f(x)+||x||$. We will show that $\overline{B}_\kappa^{||\cdot||}(0)\subseteq \overline{B}_1^{||\cdot||_x}(0)\subseteq \overline{B}_\lambda^{||\cdot||}(0)$.

Let $y\in \overline{B}_\kappa^{||\cdot||}(0)$. Hence $x+y\in \overline{B}_\kappa^{||\cdot||}(x)\subseteq K_{f,\epsilon}$. From $-y\in \overline{B}_\kappa^{||\cdot||}(0)$ we furthermore obtain $x-y\in \overline{B}_\kappa^{||\cdot||}(x)\subseteq K_{f,\epsilon}$. Thus $0\leq x+y$ and $0\leq x-y$, i.e.\ $y\in [-x,x]=\overline{B}_1^{||\cdot||_x}(0)$. 

To show the second inclusion let $y\in \overline{B}_1^{||\cdot||_x}(0)=[-x,x]$. The linear functional $f$ is positive as presented in Proposition \ref{pro:basicpropertiesoftheicecreamcone}. Thus $f(y)\leq f(x)$. From $x+y\in K_{f,\epsilon}$ we obtain $\epsilon ||x+y||\leq f(x+y)=f(x)+f(y)\leq 2 f(x)$. Hence $||y||\leq ||x+y||+||x||\leq \frac{2}{\epsilon}f(x)+||x||=\lambda$, i.e.\ $y\in \overline{B}_\lambda^{||\cdot||}(0)$. 	
\end{proof}
		
\section{Decreasing nets and the order bound topology}

In this section let $X$ be a partially ordered vector space. Remark \ref{rem:basicorderconvergencestuff}(b) states that every net $(x_\alpha)_{\alpha \in A}$ in $X$ with $x_\alpha \downarrow 0$ satisfies $x_\alpha \xrightarrow{\tau_o}0$. In this section we will characterise the partially ordered vector spaces (with order units) in which $x_\alpha \downarrow 0$ implies $x_\alpha \xrightarrow{\tau_{ob}}0$. 
		
\begin{lemma} \label{lem:XsatisfiesdecreasinfTobpropertyimplieseachorderunittobenetcatching}
	Let $X$ be such that for all nets $(x_\alpha)_{\alpha \in A}$ with $x_\alpha \downarrow 0$ there is $x_\alpha \xrightarrow{\tau_{ob}}0$. Then each order unit $u$ is net catching in $X$.
\end{lemma}

\begin{proof}
	Let $u$ be an order unit of $X$ and $(x_\alpha)_{\alpha \in A}$ a net in $X$ such that $x_\alpha \downarrow 0$. Hence $[-u,u] \in \mathcal{B}_{ob}$ by Proposition \ref{pro:orderunitsgiveneighbourhoodsof0}. By the assumption we obtain $x_\alpha \xrightarrow{\tau_{ob}}0$. Since $[-u,u]$ is a $\tau_{ob}$-neighbourhood of $0$ there is $\alpha \in A$ such that $x_\alpha \in [-u,u]$, hence $x_\alpha \leq u$.
\end{proof}
		
		This allows to show the following characterization. 
		
\begin{theorem} \label{the:characterizationwhendodecreasingnetsconvergewithrespecttotheorderboundtopology}
	Let $X$ be Archimedean and contain order units. Then the following statements are equivalent.
\begin{itemize}
	\item[(i)] $X$ has net catching elements.
	\item[(ii)] $\tau_o=\tau_{ob}$.
	\item[(iii)] For all nets $(x_\alpha)_{\alpha \in A}$ with $x_\alpha \xrightarrow{\tau_o}x$ we have $x_\alpha \xrightarrow{\tau_{ob}}x$.
	\item[(iv)] For all nets $(x_\alpha)_{\alpha \in A}$ with $x_\alpha \downarrow x$ we have $x_\alpha \xrightarrow{\tau_{ob}}x$.				
\end{itemize}
\end{theorem}
		
\begin{proof} If $X$ has order units, then $X$ is directed. Thus (i) implies (ii) by Theorem \ref{the:ArchimedeananddirectedandnetcatchingthenTobisTo}. Clearly (iii) follows from (ii). If $(x_\alpha)_{\alpha \in A}$ is a net with $x_\alpha \downarrow x$, then $x_\alpha \xrightarrow{\tau_o}$ by Remark \ref{rem:basicorderconvergencestuff}(b). Hence (iv) is implied by (iii). The remaining implication can be obtained from Lemma \ref{lem:XsatisfiesdecreasinfTobpropertyimplieseachorderunittobenetcatching}. 
\end{proof}

\begin{corollary}
	Let $X$ be an infinite dimensional Archimedean vector lattice with order units. Then there is a net $(x_\alpha)_{\alpha \in A}$ in $X$ such that $x_\alpha \downarrow 0$, but not $x_\alpha \xrightarrow{\tau_{ob}}0$. Indeed, $X$ has no net catching element by Theorem \ref{the:Characterization_existencence_for_Archimedean_VL}. Thus (iv) in Theorem \ref{the:characterizationwhendodecreasingnetsconvergewithrespecttotheorderboundtopology} can not be satisfied.
\end{corollary}

\begin{acknowledgement}
The author would like to thank his supervisor Anke Kalauch for the patient guidance, encouragement and advice she has provided throughout the completion of the authors Master's thesis. A lot of the material in this paper is due to very much enjoyable discussions in the context of this thesis.  
\end{acknowledgement}

 \vspace{10mm} \noindent
\begin{tabular}{l}
Till Hauser \\
Fakult\"at f\"ur Mathematik und Informatik\\
Institut f\"ur Mathematik\\
Friedrich-Schiller-Universit\"at Jena\\
Ernst-Abbe-Platz 2\\
07743 Jena\\
Germany\\
till.hauser@uni-jena.de \end{tabular}

\end{document}